\documentclass[12pt,draft]{amsart}
\usepackage{amsmath,amssymb,amsthm,bm}
\usepackage{graphicx,enumerate}
\usepackage{layout}

\theoremstyle{plain}
\newtheorem{theorem}{Theorem}[section]
\newtheorem{lemma}[theorem]{Lemma}

\theoremstyle{definition}

\numberwithin{equation}{section}

\title[Twisted Symmetric Group Actions]{\Large Twisted Symmetric Group Actions}


\thanks{The first-named author was partially supported by Waseda University Grant
for Special Research Projects No. 2007B-067 and Rikkyo University
Special Fund for Research. The second-named author was partially
supported by National Center for Theoretic Sciences (Taipei
Office).}

\begin{document}
\maketitle
%\layout
\begin{center}
\begin{tabular}{lll}
Akinari Hoshi& &Ming-chang Kang\\
Department of Mathematics & &Department of Mathematics\\
Rikkyo University& and\hspace*{3mm} &National Taiwan University\\
Tokyo, Japan & &Taipei, Taiwan \\
E-mail: \texttt{hoshi@rikkyo.ac.jp} & &E-mail: \texttt{kang@math.ntu.edu.tw}
\end{tabular}
\end{center}

\bigskip
%%%%%%%%%%%%%%%%%%%%%%%%%%%%%%%%%%%%%%%%%%%%%%%%%%%%%%%%%%%%%%%%%%%%%%%%%%%%%%%%%%%%%%%%%%%%%
\noindent Abstract.
Let $K$ be any field, $K(x_1,\ldots,x_n)$ be the rational function field of $n$
variables over $K$, and $S_n$ and $A_n$ be the symmetric group and the alternating group
of degree $n$ respectively.
For any $a\in K\setminus\{0\}$, define an action of $S_n$ on $K(x_1,\ldots,x_n)$ by
$\sigma\cdot x_i=x_{\sigma(i)}$ for $\sigma\in A_n$ and
$\sigma\cdot x_i=a/x_{\sigma(i)}$ for $\sigma\in S_n\setminus A_n$.
Theorem. For any field $K$ and $n=3,4,5$, the fixed field $K(x_1,\ldots,x_n)^{S_n}$ is
rational (= purely transcendental) over $K$.

%%%%%%%%%%%%%%%%%%%%%%%%%%%%%%%%%%%%%%%%%%%%%%%%%%%%%%%%%%%%%%%%%%%%%%%%%%%%%%%%%%%%%%%%%%%%%%%
\section{Introduction}\label{sec-intro}

Let $K$ be any field, $K(x_1,\ldots,x_n)$ be the rational function field of $n$ variables
over $K$, and $S_n$ and $A_n$ be the symmetric group and the alternating group
of degree $n$ respectively.
For any $a\in K\setminus\{ 0\}$, define a twisted action of $S_n$ on $K(x_1,\ldots,x_n)$ by
\begin{align}
\sigma(x_i):=
\begin{cases}
\quad x_{\sigma(i)},\qquad \mathrm{if}\qquad \sigma\in A_n,\\
\quad a/x_{\sigma(i)},\quad \mathrm{if}\qquad \sigma\in S_n\setminus A_n.
\end{cases}\label{actsigma}
\end{align}

Consider the fixed subfield $K(x_1,\ldots,x_n)^{S_n}=\{\,\alpha\in K(x_1,\ldots,x_n)\,
:\, \sigma(\alpha)=\alpha\ \mathrm{for\ any}\ \sigma\in S_n\, \}$.

If $n=2$, then $K(x_1,x_2)^{S_2}=K(x_1+(a/x_2),\, ax_1/x_2)$ is rational
(= purely transcendental) over $K$.
The answer when $a=1$ (equivalently when $a\in K^{\times 2}$) is provided
by the following theorem.

\begin{theorem}[Hajja and Kang {\cite[Theorem 3.5]{HK3}}]\label{thHK}
Let $K$ be any field and $a\in K^{\times 2}$.
Then $K(x_1,\ldots,x_n)^{S_n}$ is rational over $K$.
\end{theorem}

The case when $a\in K^\times\setminus K^{\times 2}$ and $n \geq 3$
had been intangible for many years (see \cite[p.638]{HK3},
\cite[Example 5.12, p.147]{Ha} and \cite[Question 3.8,
p.215]{Ka1}); even the case $n=3$ was unknown previously. The
following Theorem 1.2 is our recent result of this question for
the cases of $n=3,4,5$.

\begin{theorem}\label{thmain}
Let $K$ be any field, $a\in K\setminus\{0\}$, and $S_n$ act on
$K(x_1,\ldots,x_n)$ as defined in Formula $(\ref{actsigma})$. If
$n=3,4,5$, then $K(x_1,\ldots,x_n)^{S_n}$ is rational over $K$.
\end{theorem}

\medskip

The proof of  Theorem \ref{thmain} will be given in Section
\ref{se2}. It is interesting to note that we use three different
methods for the three cases of $n$. It seems that there is no
unified proof for the three cases. One of the reasons to explain
this phenomenon is that the solutions to Noether's problem for the
alternating group $A_n$ are rather different when $n=3$ and when
$n=5$ (see Theorem \ref{lemMasuda} and Theorem \ref{lemMaeda}).
Since Noether's problem for $A_n$ is still an open problem when
$n\ge 6$ (see \cite{Mae} and \cite[Section 4]{HK2} for this
problem), it is not so surprising that our question is solvable at
present only for $n \leq 5$. It is still unknown whether the fixed
field $K(x_1,\ldots,x_n)^{S_n}$ of our question is rational or not
when $n\geq 6$.

In Section \ref{se3} we propose another approach to this question.
It will be shown that $K(x_1,\ldots,x_n)^{S_n}$ is isomorphic to
the function field of a conic bundle over $\mathbb{P}^{n-1}$ of
the form $x^2-ay^2=h(v_1,\ldots,v_{n-1})$ with affine coordinates
$v_1,\ldots,v_{n-1}$ (see Theorem \ref{thconic}). Although this
approach is valid only when char $K\neq 2$, it does provide a new
technique in studying rationality problems. The structure of a
conic bundle together with its rationality problem is a central
subject in algebraic geometry \cite{Is}. Fortunately, when $n=3$
and $n=4$, the conic bundle in our question contains singularities
and the rationality problem can be solved by suitable blowing-up
process. In particular, we find another proof of Theorem
\ref{thmain} when char $K\neq 2$ and $n=3,4$. For other
rationality problems of conic bundles, see \cite[Section 4]{Ka3}.

Since the fixed field $K(x_1,\ldots,x_n)^{S_n}$ is the quotient
field of the ring of invariants $K[x_1,\ldots,x_n]^{S_n}$, it
seems plausible to study the rationality problem of
$K(x_1,\ldots,x_n)^{S_n}$ through the structure of
$K[x_1,\ldots,x_n]^{S_n}$. This strategy is carried out in Section
4 and we give another proof of Theorem \ref{thmain} when char $K=
2$ and $n=3, 4$.

%%%%%%%%%%%%%%%%%%%%%%%%%%%%%%%%%%%%%%%%%%%%%%%%%%%%%%%%%%%%%%%%%%%%%%%%%%%%%%%%%%%%%%%%%%%%%%%
\section{Proof of Theorem \ref{thmain}}\label{se2}

First of all we recall two results to be used in the sequel.

\begin{theorem}[Kang {\cite[Theorem
2.4]{Ka2}}]\label{lemHK} Let $K$ be any field and $K(x,y)$ be the
rational function field of two variables over $K$. Let $\sigma$ be
a $K$-automorphism on $K(x,y)$ defined by
\begin{align*}
\sigma\,:\, x\ \mapsto\ \frac{a}{x},\quad y\ \mapsto\ \frac{b}{y}
\end{align*}
where $a\in K\setminus\{0\}$ and $b=c(x+(a/x))+d$ such that $c,d
\in K$ and at least one of $c,d$ is non-zero. Then
$K(x,y)^{\langle \sigma\rangle}=K(s,t)$ where
\begin{align*}
s=\frac{x-(a/x)}{xy-(ab/xy)},\quad t=\frac{y-(b/y)}{xy-(ab/xy)}.
\end{align*}
\end{theorem}

The following result is essentially due to Masuda \cite[p.62]{Mas}
when char $K\neq 3$ (with a misprint in the original expression).
We thank Y. Rikuna who pointed out that when char $K=3$, the same
formula is still valid, if we compare this formula with the proof
in \cite{Ku}. For the convenience of the reader, we provide a new
proof of it.

\begin{theorem}[Masuda {\cite[Theorem 3]{Mas}}]\label{lemMasuda}
Let $K$ be any field, $K(x_1,x_2,x_3)$ be the rational function
field of three variables over $K$.
Let $\sigma$ be a $K$-automorphism on $K(x_1,x_2,x_3)$ defined by
\begin{align*}
\sigma\ :\ x_1\ \longmapsto\ x_2\ \longmapsto\ x_3\ \longmapsto\ x_1.
\end{align*}
Then $K(x_1,x_2,x_3)^{\langle\sigma\rangle}=K(s_1,u,v)=K(s_3,u,v)$
where $s_i$ is the elementary symmetric function of degree $i$ for
$1 \leq i \leq 3$, and $u$ and $v$ are defined by
\begin{align*}
u&:=\frac{x_1x_2^2+x_2x_3^2+x_3x_1^2-3x_1x_2x_3}{x_1^2+x_2^2+x_3^2-x_1x_2-x_2x_3-x_3x_1},\\
v&:=\frac{x_1^2x_2+x_2^2x_3+x_3^2x_1-3x_1x_2x_3}{x_1^2+x_2^2+x_3^2-x_1x_2-x_2x_3-x_3x_1}.
\end{align*}
Moreover, we have the following identities
\begin{align*}
&s_2=s_1(u+v)-3(u^2-uv+v^2),\\
&s_3=s_1uv-(u^3+v^3),\\
&x_1x_2^2+x_2x_3^2+x_3x_1^2=s_1^2u-3s_1u^2+3(2u-v)(u^2-uv+v^2),\\
&x_1^2x_2+x_2^2x_3+x_3^2x_1=s_1^2v-3s_1v^2-3(u-2v)(u^2-uv+v^2).
\end{align*}
\end{theorem}
\begin{proof}
With the aid of computer packages, e.g. Mathematica or Maple, it
is easy to verify the formulae of $s_2$, $s_3$,
$x_1x_2^2+x_2x_3^2+x_3x_1^2$, $x_1^2x_2+x_2^2x_3+x_3^2x_1$ in
terms of $s_1$, $u$, $v$ as given in the statement of the above
theorem. We have $[K(x_1,x_2,x_3) : K(s_1,s_2,s_3)]=6$ and
$[K(x_1,x_2,x_3)^{\langle\sigma\rangle} : K(s_1,s_2,s_3)]=2$.
Since $x_1x_2^2+x_2x_3^2+x_3x_1^2 \not\in K(s_1,s_2,s_3)$, it
follows that $K(x_1,x_2,x_3)^{\langle\sigma\rangle}
=K(s_1,s_2,s_3,x_1x_2^2+x_2x_3^2+x_3x_1^2)\subset K(s_1,u,v)$.
Hence
$K(x_1,x_2,x_3)^{\langle\sigma\rangle}=K(s_1,u,v)=K(s_3,u,v)$.
\end{proof}

\bigskip

%%%%%%%%%%%%%%%%%%%%%%%%%%%%%%%%%%%%%%%%%%%%%%%%%%%%%%%%%%%%%%%%%%%%%%%%%%%%%%%%%%%%%%%%%%%%%

Proof of Theorem \ref{thmain} when $n=3$
------------------------------------------

\medskip

Let $\sigma=(1,2,3)$, $\tau=(1,2)\in S_3$.

By Theorem \ref{lemMasuda}, we find that
$K(x_1,x_2,x_3)^{\langle\sigma\rangle}=K(s_3,u,v)$.

Now $\tau(x_1)=a/x_2$, $\tau(x_2)=a/x_3$, $\tau(x_3)=a/x_3$.
Note that
\begin{align*}
&\tau(s_1)=as_2/s_3,\quad \tau(s_2)=a^2s_1/s_3,\quad \tau(s_3)=a^3/s_3,\\
&\tau(x_1x_2^2+x_2x_3^2+x_3x_1^2)=a^3(x_1x_2^2+x_2x_3^2+x_3x_1^2)/s_3^2,\\
&\tau(x_1^2x_2+x_2^2x_3+x_3^2x_1)=a^3(x_1^2x_2+x_2^2x_3+x_3^2x_1)/s_3^2.
\end{align*}
With the aid of Theorem \ref{lemMasuda}, it is not difficult to
find that
\begin{align}
\tau\ :\ &\ s_3\ \longmapsto\ \frac{a^3}{s_3},\quad
u\ \longmapsto\ \frac{au}{u^2-uv+v^2},\quad
v\ \longmapsto\ \frac{av}{u^2-uv+v^2}.\label{acttM}
\end{align}

Define $w:=u/v$.
Then $K(s_3,u,v)=K(s_3,v,w)$ and
\begin{align*}
\tau\ :\ &\ s_3\ \longmapsto\ \frac{a^3}{s_3},\quad
v\ \longmapsto\ \frac{a}{v(1-w+w^2)},\quad w\ \longmapsto\ w.
\end{align*}

By Theorem \ref{lemHK}, $K(s_3,v,w)^{\langle\tau\rangle}$ is
rational over $K(w)$. Hence $K(x_1,x_2,x_3)^{S_3} =
K(s_3,v,w)^{\langle\tau\rangle}$ is rational over $K$.\qed

\bigskip

%%%%%%%%%%%%%%%%%%%%%%%%%%%%%%%%%%%%%%%%%%%%%%%%%%%%%%%%%%%%%%%%%%%%%%%%%%%%%%%%%%%%%%%%%

Proof of Theorem \ref{thmain} when $n=4$
------------------------------------------

\medskip

Define
\begin{align*}
\sigma&:=(123)\hspace*{5mm}\ :\ x_1\longmapsto x_2\longmapsto x_3\longmapsto x_1,\\
\tau&:=(12)\hspace*{7mm}\ :\ x_1\longmapsto a/x_2,\quad x_2\longmapsto a/x_1,\
x_3\longmapsto a/x_3,\quad x_4\longmapsto a/x_4,\\
\rho_1&:=(12)(34)\ :\ x_1\longmapsto x_2,\quad x_2\longmapsto x_1,\quad
x_3\longmapsto x_4,\quad x_4\longmapsto x_3,\\
\rho_2&:=(13)(24)\ :\ x_1\longmapsto x_3,\quad x_3\longmapsto x_1,\quad
x_2\longmapsto x_4,\quad x_4\longmapsto x_2.
\end{align*}

Note that the following is a normal series:
\[
\{1\}
\ \lhd\ V_4=\langle\rho_1,\rho_2\rangle
\ \lhd\ A_4=\langle\sigma,\rho_1,\rho_2\rangle
\ \lhd\ S_4=\langle\sigma,\tau,\rho_1,\rho_2\rangle.
\]

First we will show that $K(x_1,\ldots,x_4)^{V_4}$ is rational over
$K$.

Define
\begin{align*}
s_1&:=x_1+x_2+x_3+x_4,\quad s_4:=x_1x_2x_3x_4,\\
S&:=\frac{x_1+x_2-x_3-x_4}{x_1x_2-x_3x_4},\quad
T:=\frac{x_1-x_2-x_3+x_4}{x_1x_4-x_2x_3},\quad
U:=\frac{x_1-x_2+x_3-x_4}{x_1x_3-x_2x_4}.
\end{align*}
Then we have $K(s_1,s_4,S,T,U)\subset K(x_1,x_2,x_3,x_4)^{V_4}$ and
\begin{align}
\sigma : s_1\longmapsto s_1,\quad s_4\longmapsto s_4,\quad S\longmapsto T,\quad
T\longmapsto U,\quad U\longmapsto S.\label{acts3}
\end{align}

\begin{lemma}
$K(x_1,x_2,x_3,x_4)^{V_4}=K(s_1,S,T,U)=K(s_4,S,T,U)$.
\end{lemma}
\begin{proof}
Define
\begin{align*}
u_1:=S+T+U,\quad u_2:=ST+TU+SU,\quad u_3:=STU.
\end{align*}
Then it can be checked that $K(x_1,x_2,x_3,x_4)=K(s_1,S,T,U)(x_4)$ directly from
the following equalities:
\begin{align*}
x_1=\frac{4-s_1T+(-2u_1+s_1T(S+U))x_4+SU(1-s_1T)x_4^2+u_3x_4^3}{S-T+U-SUx_4},\\
x_2=\frac{4-s_1U+(-2u_1+s_1U(T+S))x_4+TS(1-s_1U)x_4^2+u_3x_4^3}{T-U+S-TSx_4},\\
x_3=\frac{4-s_1S+(-2u_1+s_1S(U+T))x_4+UT(1-s_1S)x_4^2+u_3x_4^3}{U-S+T-UTx_4}.
\end{align*}
We see that $[K(s_1,S,T,U)(x_4):K(s_1,S,T,U)]\leq 4$ by the equality
\begin{align*}
u_1^2-4u_2+s_1u_3+(8-s_1u_1)u_3x_4-(2u_1-s_1u_2)u_3x_4^2-s_1u_3^2x_4^3+u_3^2x_4^4=0.
\end{align*}
Hence we get $K(x_1,x_2,x_3,x_4)^{V_4}=K(s_1,S,T,U)$.
It also follows from the equality $s_4=(u_1^2-4u_2+u_3s_1)/u_3^2$ that
$K(s_1,S,T,U)=K(s_4,S,T,U)$.
\end{proof}

Now we have
\[
K(x_1,x_2,x_3,x_4)^{S_4}=(K(x_1,x_2,x_3,x_4)^{V_4})^{S_4/V_4}
=K(s_4,S,T,U)^{\langle \sigma,\tau\rangle}.
\]
The action of $\langle \sigma,\tau\rangle$ on $K(s_4,S,T,U)$ is
given by
\begin{align*}
\sigma\ :\ &\ s_4\longmapsto s_4,\quad S\longmapsto T,\quad T\longmapsto U,\ U\longmapsto S,\\
\tau\ :\ &\ s_4\longmapsto \frac{a^4}{s_4},\quad S\longmapsto \frac{-S+T+U}{aTU},\quad
T\longmapsto \frac{S+T-U}{aST},\quad U\longmapsto \frac{S-T+U}{aSU}.
\end{align*}

Define
\begin{align*}
N:=\begin{cases}\displaystyle{\frac{s_4+a^2}{s_4-a^2}},\ \mathrm{if}\ \mathrm{char}\ K\neq 2,\\
\displaystyle{\frac{s_4}{s_4+a^2}},\ \mathrm{if}\ \mathrm{char}\ K=2.
\end{cases}
\end{align*}
Then we get $K(s_4,S,T,U)=K(N,S,T,U)$, $\sigma(N)=N$ and
\begin{align*}
\tau(N)=\begin{cases}-N,\quad\ \mathrm{if}\ \mathrm{char}\ K\neq 2,\\
N+1,\ \mathrm{if}\ \mathrm{char}\ K=2.
\end{cases}
\end{align*}
Applying \cite[Theorem 1]{HK2}, we find that
$K(x_1,x_2,x_3,x_4)^{S_4}=K(N,S,T,U)^{\langle \sigma,\tau\rangle}$
is rational over $K$, provided that $K(S,T,U)^{\langle
\sigma,\tau\rangle}$ is rational over $K$. Explicitly, define $P$
by

\begin{align*}
P:=\begin{cases}
\displaystyle{N\cdot\Bigl(S+T+U+\frac{S^2+T^2+U^2-2(ST+TU+US)}{aSTU}\Bigr)},
\ \mathrm{if}\ \mathrm{char}\ K\neq 2,\\
\displaystyle{N+\frac{S+T+U}{S+T+U+aSTU}},\ \mathrm{if}\ \mathrm{char}\ K=2
\end{cases}
\end{align*}
then $K(N,S,T,U)=K(P,S,T,U)$ and
$K(x_1,x_2,x_3,x_4)^{S_4}=K(P,S,T,U)^{\langle \sigma,\tau\rangle}
=K(S,T,U)^{\langle \sigma,\tau\rangle}(P)$ where $\sigma (P) =
\tau(P)=P$.

Thus it remains to prove the following theorem.

\begin{theorem}
Let $K$ be any field and $K(S,T,U)$ be the rational function field of three variables $S,T,U$
over $K$.
Let $\sigma$ and $\tau$ be $K$-automorphisms of $K(S,T,U)$ defined by
\begin{align*}
\sigma\ :\ &\ S\longmapsto T,\quad T\longmapsto U,\quad U\longmapsto S,\\
\tau\ :\ &\ S\longmapsto \frac{-S+T+U}{aTU},\quad
T\longmapsto \frac{S+T-U}{aST},\quad U\longmapsto \frac{S-T+U}{aSU}
\end{align*}
where $a\in K\setminus\{0\}$.
Then $\langle\sigma,\tau\rangle\cong S_3$ and
$K(S,T,U)^{\langle\sigma,\tau\rangle}$ is rational over $K$.
\end{theorem}
\begin{proof}
By Theorem \ref{lemMasuda}, we may choose a transcendence basis of
$K(S,T,U)^{\langle\sigma\rangle}$ over $K$ by
\[
K(S,T,U)^{\langle\sigma\rangle}=K(f,g,h)
\]
where
\begin{align*}
f\ &=\ S+T+U,\\
g\ &=\ \frac{ST^2 + TU^2 + US^2 - 3STU}{S^2  + T^2 + U^2 - ST - TU - US},\\
h\ &=\ \frac{S^2T + T^2U + U^2S - 3STU}{S^2  + T^2 + U^2 - ST  - TU - US}.
\end{align*}
Thus we have $K(S,T,U)^{\langle\sigma,\tau\rangle}=
(K(S,T,U)^{\langle\sigma\rangle})^{\langle\tau\rangle}=K(f,g,h)^{\langle\tau\rangle}$.
The action of $\tau$ on $K(f,g,h)$ is given by
\begin{align*}
\tau\ :\ &f\longmapsto \frac{f^2-4f(g+h)+12X}{aY},\\
&g\longmapsto \frac{-f^2h(f-4h)+2f(f-2g-8h)X+24X^2-8gY}{a(f^2-2f(g+h)+4X)Y},\\
&h\longmapsto \frac{-f^2(fg+4h^2)+6f(f-2g)X+24X^2-4(f+2h)Y}{a(f^2-2f(g+h)+4X)Y}
\end{align*}
where $X=g^2-gh+h^2$ and $Y=g^3-fgh+h^3$.\\

Case 1. char $K\neq 2$. \\

Define
\[
F:=g+h,\quad \ G:=g-h,\quad \ H:=f-(g+h).
\]
Then $K(S,T,U)^{\langle\sigma\rangle}=K(f,g,h)=K(F,G,H)$ and $\tau$ acts on $K(F,G,H)$ by
\begin{align*}
&F\ \longmapsto\ \frac{4(27G^4-7FG^2H+5G^2H^2-FH^3)}{a(4FG^2-F^2H+G^2H)(3G^2+H^2)},\\
&G\ \longmapsto\ \frac{4G(FG^2+7G^2H-FH^2+H^3)}{a(4FG^2-F^2H+G^2H)(3G^2+H^2)},\\
&H\ \longmapsto\ \frac{4H(FG^2+7G^2H-FH^2+H^3)}{a(4FG^2-F^2H+G^2H)(3G^2+H^2)}.
\end{align*}
Note that $\tau(G/H)=G/H$.
Define
\[
A:= F/G,\quad B:= G,\quad C:=G/H.
\]
Then $K(S,T,U)^{\langle\sigma\rangle}=K(F,G,H)=K(A,B,C)$ and $\tau$ acts on $K(A,B,C)$ by
\begin{align*}
&A\longmapsto \frac{-A+5C-7AC^2+27C^3}{1-AC+7C^2+AC^3},\\
&B\longmapsto \frac{4(1-AC+7C^2+AC^3)}{aB(1-A^2+4AC)(1+3C^2)},\quad
C\longmapsto C.
\end{align*}
Define
\[
D:=1-AC+7C^2+AC^3,\quad E:=2C(C^2-1)/B.
\]
Then $K(A,B,C)=K(C,D,E)$ and the action of $\tau$ on $K(C,D,E)$ is given by
\begin{align*}
&C\longmapsto C,\quad D\longmapsto \frac{(1+3C^2)^3}{D},\\
&E\longmapsto -a(1+3C^2)\Bigl(D+\frac{(1+3C^2)^3}{D}-2(1+5C^2+2C^4)\Bigr)\Big/E.
\end{align*}
Hence the assertion follows from Theorem \ref{lemHK}.\\

Case 2. char $K=2$.\\

The action of $\tau$ on $K(f,g,h)$ is given by
\begin{align*}
\tau\ :\ &f\longmapsto \frac{f^2}{aY},\quad
g\longmapsto \frac{fh}{aY},\quad
h\longmapsto \frac{fg}{aY}
\end{align*}
where $Y=g^3+fgh+h^3$.
Define
\[
A:= f/(g+h),\quad B:=g/h,\quad C:=1/h.
\]
Then $K(f,g,h)=K(A,B,C)$ and $\tau$ acts on $K(A,B,C)$ by
\begin{align*}
&A\longmapsto A,\quad
B\longmapsto \frac{1}{B},\quad
C\longmapsto \frac{a}{A}\Bigl(B+\frac{1}{B}+A+1\Bigr)\Big/C.
\end{align*}
Hence the assertion follows from Theorem \ref{lemHK}.

We will give another proof when $n=4$ and char $K=2$ in Section 4.
\end{proof}

\bigskip

%%%%%%%%%%%%%%%%%%%%%%%%%%%%%%%%%%%%%%%%%%%%%%%%%%%%%%%%%%%%%%%%%%%%%%%%%%%%%%%%%%%%%%%%%%%%%

Proof of Theorem \ref{thmain} when $n=5$
------------------------------------------

\medskip

We recall Maeda's Theorem for $A_5$ action \cite{Mae}.

\begin{theorem}[Maeda \cite{Mae}]\label{lemMaeda}
Let $K$ be any field, $K(x_1,\ldots,x_5)$ be the rational function
field of five variables over $K$. Then $K(x_1,\ldots,x_5)^{A_5}$
is rational over $K$. Moreover a transcendental basis
$F_1,\ldots,F_5$ of $K(x_1,\ldots,x_5)^{A_5}$ over $K$
may be given explicitly as follows:\\
{\rm (i)} When char $K\neq 2$;
\begin{align*}
F_1&=\frac{\sum_{\sigma\in S_5}\sigma([12][13][14][15][23]^4[45]^4x_1)}
{\sum_{\sigma\in S_5}\sigma([12][13][14][15][23]^4[45]^4)},\\
F_2&=\frac{\sum_{\sigma\in S_5}\sigma([12]^3[13]^3[14]^3[15]^3[23]^{10}[45]^{10})}
{\prod_{i<j}[ij]^2\cdot\sum_{\sigma\in S_5}\sigma([12][13][14][15][23]^4[45]^4)},\\
F_3&=\frac{\sum_{\sigma\in S_5}\sigma([12]^3[13]^3[14]^3[15]^3[23]^{10}[45]^{10}x_1)}
{\prod_{i<j}[ij]^2\cdot\sum_{\sigma\in S_5}\sigma([12][13][14][15][23]^4[45]^4)},\\
F_4&=\frac{\sum_{\mu\in R_1}\mu([12]^2[13]^2[23]^2[45]^4)}{\prod_{i<j}[ij]},\\
F_5&=\frac{\sum_{\mu\in R_1}\mu([12]^2[13]^2[23]^2[14]^4[24]^4[34]^4[15]^4[25]^4[35]^4)}
{\prod_{i<j}[ij]^3}
\end{align*}
where $[ij]=x_i-x_j$ and
$R_1=\{1,$ $(34)$, $(354)$, $(234)$, $(2354)$, $(24)(35)$, $(1234)$, $(12354)$,
$(124)(35)$, $(13524)\}$.

{\rm (ii)} When char $K=2$;
\begin{align*}
F_1&=\frac{\sum_{i<j<k}x_ix_jx_k}{\sum_{i<j}x_ix_j},\quad
F_2=\frac{\sum_{i=1}^5([12][13][14][15]\cdot I^2)^{(1i)}}
{\prod_{i<j}[ij]\cdot\sum_{i<j}x_ix_j},\\
F_3&=\frac{\sum_{i=1}^5([12][13][14][15]\cdot I^2\cdot x_1)^{(1i)}}
{\prod_{i<j}[ij]\cdot\sum_{i<j}x_ix_j},\\
F_4&=\frac{\sum_{\nu\in R_3}\nu([12]^2[34]^2[13][24][15][25][35][45])}{\prod_{i<j}[ij]},\quad
F_5={\rm the\ same\ } F_5 {\rm\ as\ in\ (i)}
\end{align*}
where $[ij]=x_i-x_j$, $I=\sum_{\tau\in R_2}\tau(x_2x_3(x_2x_3+x_4^2+x_5^2))$,
$R_2=\{1$, $(34)$, $(354)$, $(234)$, $(2354)$, $(24)(35)\}$ and
$R_3=\{1$, $(234)$, $(243)$, $(152)$, $(15234)$, $(15243)$, $(125)$, $(12345)$, $(12435)$,
$(15432)$, $(154)$, $(15423)$, $(15342)$, $(15324)$, $(153)\}$.
\end{theorem}

In the above theorem, note that $R_1$, $R_2$ and $R_3$ are coset
representatives with respect to various subgroups:

\[
S_5=\bigcup_{\mu\in R_1} H_1\mu,\quad H=\bigcup_{\tau\in R_2}
H_2\tau,\quad A_5=\bigcup_{\nu\in R_3} H_3\nu
\]
where $H_1=\langle(12),(13),(45)\rangle$ $\cong$ $D_6$,
$H=\langle(23),(24),(25)\rangle$ $\cong$ $S_4$,
$H_2=\langle(23),(45)\rangle$ $\cong$ $V_4$,
$H_3=\langle(12)(34),(13)(24)\rangle$ $\cong$ $V_4$, and
$D_6$ is the dihedral group of order $12$.

\bigskip

Now we start to prove Theorem \ref{thmain} when $n=5$.

Let $\tau=(12)\in S_5$.

By Theorem \ref{lemMaeda}, we see that
$K(x_1,\ldots,x_5)^{A_5}=K(F_1,\ldots,F_5)$.

With the aid of computer packages, e.g. Mathematica or Maple,
we can evaluate the action of $\tau$ on $K(F_1,\ldots,F_5)$ as follows:
\begin{align*}
\tau\ :\
F_1\ &\longmapsto\ \frac{a}{F_1},\quad
F_2\ \longmapsto\ \frac{F_3}{F_1},\quad
F_3\ \longmapsto\ \frac{aF_2}{F_1},\\
F_4\ &\longmapsto\ -F_4,\quad
F_5\ \longmapsto\ -F_5,\ {\rm when\ char}\ K\neq 2;\\
%%%%%%%%%%%%%%%%%%%%%%%%%%%%%%%%%%%%%%%%%%%%%%%%%%%%%%%%%%
\tau\ :\
F_1\ &\longmapsto\ \frac{a}{F_1},\quad
F_2\ \longmapsto\ \frac{F_3}{F_1},\quad
F_3\ \longmapsto\ \frac{aF_2}{F_1},\\
F_4\ &\longmapsto\ F_4+1,\quad
F_5\ \longmapsto\ F_5,\ {\rm when\ char}\ K=2.
\end{align*}

Case 1. char $K\neq 2$.

Define
\[
G_1:=F_1,\ G_2:=\frac{F_4+1}{F_4-1},\
G_3:=F_4\Bigl(F_2-\frac{F_3}{F_1}\Bigr),\
G_4:=F_2+\frac{F_3}{F_1},\ G_5:=F_4F_5.
\]
Then we have $K(x_1,\ldots,x_5)^{A_5}=K(F_1,\ldots,F_5)=K(G_1,\ldots,G_5)$ and
\begin{align*}
\tau\ :\
G_1\ &\longmapsto\ \frac{a}{G_1},\quad
G_2\ \longmapsto\ \frac{1}{G_2},\quad
G_3\ \longmapsto\ G_3,\quad
G_4\ &\longmapsto\ G_4,\quad
G_5\ \longmapsto\ G_5.
\end{align*}
Hence it follows from Theorem \ref{lemHK} that
$K(x_1,\ldots,x_5)^{S_5}=K(G_3,G_4,G_5)(G_1,G_2)^{\langle\tau\rangle}$
is rational over $K$.

\medskip

Case 2. char $K=2$.

Define
\[
G_1:=F_1,\quad G_2:=F_2,\quad
G_3:=\frac{F_2F_3}{F_1},\quad
G_4:=F_4+\frac{F_3}{F_1F_2+F_3},\quad G_5:=F_5.
\]
Then we have $K(x_1,\ldots,x_5)^{A_5}=K(F_1,\ldots,F_5)=K(G_1,\ldots,G_5)$ and
\begin{align*}
\tau\ :\
G_1\ &\longmapsto\ \frac{a}{G_1},\quad
G_2\ \longmapsto\ \frac{G_3}{G_2},\quad
G_3\ \longmapsto\ G_3,\quad
G_4\ &\longmapsto\ G_4,\quad
G_5\ \longmapsto\ G_5.
\end{align*}
We use Theorem \ref{lemHK} and find that
$K(x_1,\ldots,x_5)^{S_5}=K(G_3,G_4,G_5)(G_1,G_2)^{\langle\tau\rangle}$
is rational over $K$. \qed

\bigskip

%%%%%%%%%%%%%%%%%%%%%%%%%%%%%%%%%%%%%%%%%%%%%%%%%%%%%%%%%%%%%%%%%%%%%%%%%%%%%%%%%%%%%%%%%%%%%%%

\section{Conic bundles : another approach when char $K\neq 2$}\label{se3}

Throughout this section we assume that char $K\neq 2$.

In this section, we will give another proof of Theorem
\ref{thmain} when $n=3,4$ (and char $K\neq 2$) by presenting
$K(x_1,\ldots,x_n)^{S_n}$ as the function field of a conic bundle
over $\mathbb{P}^{n-1}$.

Consider the action of $S_n$ on $K(x_1,\ldots,x_n)$ defined by Formula (\ref{actsigma}).
Because of Theorem \ref{thHK}, we may assume that $a\in K^\times\setminus K^{\times 2}$
without loss of generality.

Define $\alpha:=\sqrt{a}$ and $\mathrm{Gal}(K(\alpha)/K)=\langle\rho\rangle$ where
$\rho(\alpha)=-\alpha$.
Extend the actions of $S_n$ and $\rho$ to $K(\alpha)(x_1,\ldots,x_n)=K(\alpha)\otimes_K
K(x_1,\ldots,x_n)$ by requiring that $S_n$ acts trivially on $K(\alpha)$ and $\tau$ acts
trivially on $K(x_1,\ldots,x_n)$.

Define $z_i:=(\alpha-x_i)/(\alpha+x_i)$ for $1\leq i\leq n$.
We find that $K(\alpha)(x_1,\ldots,x_n)=K(\alpha)(z_1,\ldots,z_n)$ and
\begin{align*}
\sigma\ :\ z_i\ \longmapsto\ -z_{\sigma(i)}
\end{align*}
for any $\sigma\in S_n\setminus A_n$, and
\begin{align*}
\rho\ :\ \alpha\ \longmapsto\ -\alpha,\quad z_i\ \longmapsto\ 1/z_i.
\end{align*}

Define $z_0:=z_1+\cdots+z_n$, $y_i:=z_i/z_0$ for $1\leq i\leq n$.
Hence $y_1+\cdots+y_n=1$.

Let $t_1,\ldots,t_n$ be the elementary symmetric functions of
$y_1,\ldots,y_n$. In particular, $t_1=1$. Define
$\Delta:=\prod_{1\leq i<j\leq n}(y_i-y_j)\in K(y_1,\ldots,y_n)$
and $u:=z_0\cdot \Delta$. Note that $\Delta^2$ can be written as a
polynomial in $t_1,\ldots,t_n$, and thus in $t_2,\ldots,t_n$.

\begin{lemma}\label{lemactt}
$K(x_1,\ldots,x_n)^{S_n}=K(\alpha)(t_2,\ldots,t_n,u)^{\langle\rho\rangle}$ and
\begin{align*}
\rho\ :\ \alpha\ \longmapsto\ -\alpha,\quad
t_i \longmapsto\ t_{n-i}\, (t_n/t_{n-1})^i\, t_n^{-1},\quad
u\ \longmapsto\ f(t_2,\ldots,t_n)\cdot u^{-1}
\end{align*}
where $f(t_2,\ldots,t_n)\in K(t_2,\ldots,t_n)$ is given by
\begin{align}
f(t_2,\ldots,t_n):=(-1)^{n(n-1)/2}\,
t_n^{-(n-1)}\, (t_n/t_{n-1})^{(n+1)(n-2)/2}\,\Delta^2\label{eqf}
\end{align}
and we adopt the convention that $t_0=t_1=1$.
\end{lemma}

\begin{proof}

Note that $K(\alpha)(y_1,\ldots,y_n,z_0)
=K(\alpha)(y_1,\ldots,y_n,u)$. Since $u$ is fixed by the action of
$S_n$, it follows that $K(\alpha)(y_1,\ldots,y_n,z_0)^{S_n}
=K(\alpha)(y_1,\ldots,y_n)^{S_n}(u) =K(\alpha)(t_2,\ldots,t_n,u)$
(the last equality follows, for example, from the proof of
\cite[Lemma 1]{HK2} because $\sigma(y_i)=y_{\sigma(i)}$ for any
$\sigma\in S_n$, any $1\leq i\leq n$ ).

Thus
$K(x_1,\ldots,x_n)^{S_n}=\{K(\alpha)^{\langle\rho\rangle}(x_1,\ldots,x_n)\}^{S_n}
=K(\alpha)(x_1,\ldots,x_n)^{\langle S_n,\, \rho\rangle}
=\{K(\alpha)(x_1,\ldots,x_n)^{ S_n} \}^{\langle \rho\rangle}
=K(\alpha)(t_2,\ldots,t_n,u)^{\langle\rho\rangle}$.

It is not difficult to verify the action of $\rho$ on
$K(\alpha)(t_2,\ldots,t_n,u)$ is given as the form in the
statement of this lemma.
\end{proof}
We write
\[
\begin{cases}
\ n=2m+1,\quad \mathrm{if}\ n\ \mathrm{is\ odd},\\
\ n=2m,\hspace*{11mm}\mathrm{if}\ n\ \mathrm{is\ even}.
\end{cases}
\]
Define
\begin{align}
u_i:=t_{i+1},\quad u_{n-i}:=\rho(t_{i+1})\,
=\,t_{n-(i+1)}\,t_n^i/t_{n-1}^{i+1},\quad (i=1,\ldots,m-1)\label{defu1}
\end{align}
and
\begin{align}
\begin{cases}
\ u_m:=t_{m+1},\quad u_{m+1}:=\rho(t_{m+1})=t_m\,t_n^m/t_{n-1}^{m+1},
\quad \mathrm{if}\ n\ \mathrm{is\ odd},\\
\ u_m:=t_n/t_{n-1},
\hspace*{62mm} \mathrm{if}\ n\ \mathrm{is\ even}.
\end{cases}\label{defu2}
\end{align}\\
Then we obtain

\begin{lemma}\label{lemactu}
$K(x_1,\ldots,x_n)^{S_n}=K(\alpha)(u_1,\ldots,u_{n-1},u)^{\langle\rho\rangle}$
and
\begin{align*}
\rho\ :\ &\alpha\ \longmapsto\ -\alpha,\quad u_i\ \longmapsto\ u_{n-i},
\quad (i=1,\ldots,n-1),\\
\ \ &u\ \longmapsto\ g(u_1,\ldots,u_{n-1})\cdot u^{-1}
\end{align*}
where $g(u_1,\ldots,u_{n-1})=f(t_2,\ldots,t_n)$ and $f(t_2,\ldots,t_n)$ is given as
in $(\ref{eqf})$.
\end{lemma}

\begin{proof}
The assertion follows from
$K(\alpha)(t_2,\ldots,t_n,u)=K(\alpha)(u_1,\ldots,u_{n-1},u)$ and
Lemma \ref{lemactt}. Indeed we may show $K(t_2,\ldots,t_n)\subset
K(u_1,\ldots,u_{n-1})$ as follows.

Case 1. $n=2m+1$ is odd.

The fact that $t_2,\ldots,t_{m+1}\in K(u_1,\ldots,u_{n-1})$
follows from Formulae (3.2) and (3.3).

$t_n\in K(u_1,\ldots,u_{n-1})$ because
\begin{align*}
\Bigl(\frac{u_m^{m+1}}{u_{m-1}^m}\Bigr)\,u_{m+1}^m\,\Bigl(\frac{1}{u_{m+2}}\Bigr)^{m+1}
=\Bigl(\frac{t_{m+1}^{m+1}}{t_m^m}\Bigr) \Bigl(\frac{t_m\,t_n^m}{t_{n-1}^{m+1}}\Bigr)^m
\Bigl(\frac{t_{n-1}^m}{t_{m+1}\,t_n^{m-1}}\Bigr)^{m+1}=t_n.
\end{align*}

$t_{n-1}\in K(u_1,\ldots,u_{n-1})$ because
\begin{align*}
t_n \Bigl(\frac{u_{m-1}}{u_m}\Bigr)\,u_{m+2}\,\Bigl(\frac{1}{u_{m+1}}\Bigr)
=t_n \Bigl(\frac{t_m}{t_{m+1}}\Bigr) \Bigl(\frac{t_{m+1}\,t_n^{m-1}}{t_{n-1}^m}\Bigr)
\Bigl(\frac{t_{n-1}^{m+1}}{t_m\,t_n^m}\Bigr)=t_{n-1}.
\end{align*}

From Formula (3.2) we find that
$t_{n-(i+1)}=u_{n-i}\,t_{n-1}^{i+1}/t_n^i$ for $1\leq i\leq m-2$.
Thus $t_{m+2},\ldots,t_{n-2}\in K(u_1,\ldots,u_{n-1})$.

\medskip

Case 2. $n=2m$ is even.

$t_2,\ldots,t_m\in K(u_1,\ldots,u_{n-1})$ follows from Formula
(3.2).

From Formulae (3.2) and (3.3), we get
\begin{align*}
\frac{u_{k+1}}{u_{k+2}}=\frac{t_k}{t_{k+1}}\cdot\frac{t_n}{t_{n-1}}
=\frac{t_k}{t_{k+1}}\cdot u_m,\quad
\end{align*}
where $k=m,\ldots,2m-3$. We find that
$t_{k+1}=t_k\,u_m\,u_{k+2}/u_{k+1}\in K(u_1,\ldots,u_{n-1})$ for
$m\leq k\leq 2m-3$.

From Formula (3.2), we have $u_{n-1}=t_{n-2}\,
t_n/t_{n-1}^2=t_{n-2}\, u_m/t_{n-1}$. Hence $t_{n-1}=t_{n-2}\,
u_m/u_{n-1}\in K(u_1,\ldots,u_{n-1})$.

Since $t_n=u_m\,t_{n-1}$, it follows that $t_n\in
K(u_1,\ldots,u_{n-1})$.

\end{proof}

\bigskip

We will change the variables $u_1,\ldots,u_{n-1}$ by
$v_1,\ldots,v_{n-1}$ defined as follows.

When $n=2m+1$ is odd, define
\begin{align*}
v_i:=\frac{u_i+u_{n-i}}{2},\quad
v_{n-i}:=\frac{\alpha(u_i-u_{n-i})}{2},\quad
\end{align*}
where $i=1,\ldots,m$.

When $n=2m$ is even, define
\begin{align*}
v_m:=u_m,\quad v_i:=\frac{u_i+u_{n-i}}{2},\quad
v_{n-i}:=\frac{\alpha(u_i-u_{n-i})}{2},
\end{align*}
where $i=1,\ldots,m-1$.

Thus
$K(\alpha)(u_1,\ldots,u_{n-1},u)=K(\alpha)(v_1,\ldots,v_{n-1},u)$.

By Lemma \ref{lemactu}, we get
\begin{lemma}\label{lemactv}
$K(x_1,\ldots,x_n)^{S_n}=K(\alpha)(v_1,\ldots,v_{n-1},u)^{\langle\rho\rangle}$
and
\begin{align*}
\rho\ :\ &\alpha\ \longmapsto\ -\alpha,\quad v_i\ \longmapsto\ v_i,\
(i=1,\ldots,n-1),\quad u\ \longmapsto\ h(v_1,\ldots,v_{n-1})\cdot u^{-1}
\end{align*}
where $h(v_1,\ldots,v_{n-1})=f(t_2,\ldots,t_n)$
and $f(t_2,\ldots,t_n)$ is given as in $(\ref{eqf})$.
\end{lemma}

Hence we get the following theorem which asserts that
$K(x_1,\ldots,x_n)^{S_n}$ is the function field of a conic bundle
over $\mathbb{P}^{n-1}$ of the form
$x^2-ay^2=h(v_1,\ldots,v_{n-1})$ with affine coordinates
$v_1,\ldots,v_{n-1}$ (see, for example, \cite[p.73]{Sh} for conic
bundles over $\mathbb{P}^1$).

\begin{theorem}\label{thconic}
$K(x_1,\ldots,x_n)^{S_n}=K(x,y,v_1,\ldots,v_{n-1})$ and the generators $x$, $y$,
$v_1,\ldots,v_{n-1}$
satisfy the relation $x^2-ay^2=h(v_1,\ldots,v_{n-1})$ where
$h(v_1,\ldots,v_{n-1})=f(t_2,\ldots,t_n)$ and $f(t_2,\ldots,t_n)$
is given as in $(\ref{eqf})$.
\end{theorem}

\begin{proof}
Define
\[
x:=\frac{1}{2}\Bigl(u+\frac{h(v_1,\ldots,v_{n-1})}{u}\Bigr),\quad
y:=\frac{1}{2\alpha}\Bigl(u-\frac{h(v_1,\ldots,v_{n-1})}{u}\Bigr).
\]
Then we obtain $K(x,y,v_1,\ldots,v_{n-1})\subset K(x_1,\ldots,x_n)^{S_n}
=K(\alpha)(v_1,\ldots,v_{n-1},u)$.
Thus $K(x,y,v_1,\ldots,v_{n-1})$ $=$ $K(x_1,\ldots,x_n)^{S_n}$ follows from
$K(x,y,v_1,\ldots,v_n)(u)=K(\alpha)(v_1,\ldots,v_{n-1},u)$ and
$[K(x,y,v_1,\ldots,v_n)(u):K(x,y,v_1,\ldots,v_n)]=2$.
We also have $x^2-ay^2=h(v_1,\ldots,v_{n-1})$ by the definition.
\end{proof}

\bigskip
%%%%%%%%%%%%%%%%%%%%%%%%%%%%%%%%%%%%%%%%%%%%%%%%%%%%%%%%%%%%%%%%%%%%%%%%%%%%%%%%%%%%%%%%%%%%%
Proof of Theorem \ref{thmain} when $n=3$ and char $K\neq 2$ ------------------------------------

\medskip

Step $1$. By Lemma \ref{lemactt} we find that
$K(x_1,x_2,x_3)^{S_3}=K(\alpha)(t_2,t_3,u)^{\langle\rho\rangle}$ where
\begin{align*}
\rho\ :\ \alpha\ \longmapsto\ -\alpha,\quad t_2\ \longmapsto\ t_2^{-2}\,t_3,\quad
t_3\ \longmapsto\ t_2^{-3}\,t_3^2,\quad u\ \longmapsto\ -t_2^{-2}\,\Delta^2\cdot u^{-1}.
\end{align*}

Note that $\Delta^2=\prod_{1\leq i<j\leq 3}(y_i-y_j)^2=t_2^2-4t_2^3-4t_3+18t_2t_3-27t_3^2$
because $t_1=1$.

Define $u_1:=t_2$, $u_2:=\rho(t_2)=t_2^{-2}t_3$.
Then $K(\alpha)(t_2,t_3,u)=K(\alpha)(u_1,u_2,u)$ and
\[
\rho\ :\ u_1\ \longmapsto\ u_2\ \longmapsto\ u_1,\quad u\ \longmapsto\ g(u_1,u_2)\cdot u^{-1}
\]
where $g(u_1,u_2)=-1+4u_1+4u_2-18u_1u_2+27u_1^2u_2^2$.

Define $v_1:=(u_1+u_2)/2$, $v_2:=\alpha(u_1-u_2)/2$.
Then $\rho\ :\ v_1\ \longmapsto\ v_1,\ v_2\ \longmapsto\ v_2$ and
$g(u_1,u_2)=h(v_1,v_2)$ where
\[
h(v_1,v_2)=-1+8v_1-18v_1^2+27v_1^4+(18/a)v_2^2-(54/a)v_1^2v_2^2+(27/a^2)v_2^4.
\]

Hence $K(x_1,x_2,x_3)^{S_3}=K(\alpha)(v_1,v_2,u)^{\langle\rho\rangle}=K(x,y,v_1,v_2)$ where
\[
x=\frac{1}{2}\Bigl(u+\frac{h(v_1,v_2)}{u}\Bigr),\quad
y=\frac{1}{2\alpha}\Bigl(u-\frac{h(v_1,v_2)}{u}\Bigr).
\]
Note that $x$ and $y$ satisfy the relation
\begin{align}
x^2-ay^2&=h(v_1,v_2)\label{relv1}\\
&=(1+v_1)(-1+3v_1)^3-(18/a)v_2^2(-1+3v_1^2)+(27/a^2)v_2^4.\nonumber
\end{align}

%%%%%%%%%%%%%%%%%%%%%%%%%%%%%%%%%%%%%%%%%%%%%%%%%%%%%%%%%%%%%%%%%%%%%%%%%%%%%%%%%%%%%%%%%%%%%
Step $2$.
Suppose that char $K=3$.
Then the relation (\ref{relv1}) becomes
$x^2-ay^2=-1-v_1$.
Hence $K(x_1,x_2,x_3)^{S_3}=K(x,y,v_1,v_2)=K(x,y,v_2)$ is rational over $K$.

\bigskip
%%%%%%%%%%%%%%%%%%%%%%%%%%%%%%%%%%%%%%%%%%%%%%%%%%%%%%%%%%%%%%%%%%%%%%%%%%%%%%%%%%%%%%%%%%%%%
Step $3$.
From now on, we assume that char $K\neq 2,3$.

We normalize the generators $v_1,v_2$ by defining
\begin{align*}
T_1:=3v_1,\quad T_2:=3v_2/a.
\end{align*}
We get $K(x_1,x_2,x_3)^{S_3}=K(x,y,T_1,T_2)$ with a relation
\begin{align}
3x^2-3ay^2=-3+8T_1-6T_1^2+T_1^4+6aT_2^2-2aT_1^2T_2^2+a^2T_2^4.\label{relX}
\end{align}

\bigskip
%%%%%%%%%%%%%%%%%%%%%%%%%%%%%%%%%%%%%%%%%%%%%%%%%%%%%%%%%%%%%%%%%%%%%%%%%%%%%%%%%%%%%%%%%%%%%
Step $4$.
Find the singularities of (\ref{relX}).
We get $x=y=-1+T_1=T_2=0$.
Define $T_3:=-1+T_1$.
The relation (\ref{relX}) becomes
\begin{align}
3x^2-3ay^2=4aT_2^2+a^2T_2^4-4aT_2^2T_3-2aT_2^2T_3^2+4T_3^3+T_3^4.\label{relT}
\end{align}

\bigskip
%%%%%%%%%%%%%%%%%%%%%%%%%%%%%%%%%%%%%%%%%%%%%%%%%%%%%%%%%%%%%%%%%%%%%%%%%%%%%%%%%%%%%%%%%%%%%
Step $5$.
Blow-up the equation (\ref{relT}), i.e. define $X_2:=x/T_3$, $Y_2:=y/T_3$, $T_4:=T_2/T_3$.
Then $K(x,y,T_1,T_2)=K(x,y,T_2,T_3)=K(X_2,Y_2,T_3,T_4)$ and the relation (\ref{relT})
becomes
\begin{align}
3X_2^2-3aY_2^2&=4T_3+T_3^2+4aT_4^2-4aT_3T_4^2-2aT_3^2T_4^2+a^2T_3^2T_4^4\label{relY}\\
&=(T_3-aT_3T_4^2)^2+4(T_3-aT_3T_4^2)+4aT_4^2\nonumber\\
&=(T_3-aT_3T_4^2)(4+T_3-aT_3T_4^2)+4aT_4^2.\nonumber
\end{align}

Define
\begin{align*}
X_3&:=\frac{X_2}{T_3-aT_3T_4^2},\quad Y_3:=\frac{Y_2}{T_3-aT_3T_4^2},\\
S_1&:=\frac{4+T_3-aT_3T_4^2}{T_3-aT_3T_4^2},\quad S_2:=\frac{T_4}{T_3-aT_3T_4^2}.
\end{align*}

Note that $K(X_2,Y_2,T_3,T_4)=K(X_3,Y_3,S_1,S_2)$.
For $S_1\in K(X_3,Y_3,S_1,S_2)$ and $S_1$ is a fractional linear transformation
of the ``variable'' $T_3-aT_3T_4^2$.
Hence $T_3-aT_3T_4^2\in K(X_3,Y_3,S_1,S_2)$.
Thus $T_4=S_2\cdot(T_3-aT_3T_4^2)\in K(X_3,Y_3,S_1,S_2)$ also.
Now $S_1$ is a fractional linear transformation of $T_3$ with coefficients in $K(T_4)$.
Hence $T_3\in K(X_3,Y_3,S_1,S_2)$.
It follows that $X_2,Y_2\in K(X_3,Y_3,S_1,S_2)$ also.

The relation (\ref{relY}) becomes
\begin{align*}
3X_3^2-3aY_3^2=S_1+4aS_2^2
\end{align*}
which is linear in $S_1$.
Hence $K(x_1,x_2,x_3)^{S_3}=K(X_3,Y_3,S_1,S_2)=K(X_3,Y_3,S_2)$ is rational over $K$.

\bigskip
%%%%%%%%%%%%%%%%%%%%%%%%%%%%%%%%%%%%%%%%%%%%%%%%%%%%%%%%%%%%%%%%%%%%%%%%%%%%%%%%%%%%%%%%%%%%%
Step $6$. Here is another proof.
Instead of the method in Step $5$, we may proceed as follows:

Define $X_4:=x/T_3^2$, $Y_4:=y/T_3^2$, $T_4:=T_2/T_3$, $T_5:=1/T_3$.
Then $K(x,y,T_2,T_3)$ $=K(X_4,Y_4,T_4,T_5)$ and the relation
(\ref{relT}) becomes
\begin{align}
3X_4^2-3aY_4^2=1-2aT_4^2+a^2T_4^4+4T_5-4aT_4^2T_5+4aT_4^2T_5^2.\label{relYT}
\end{align}

The singularities of (\ref{relYT}) are $X_4=Y_4=T_4\pm (1/\sqrt{a})=T_5=0$.
Blow-up with respect to $1-aT_4^2$, i.e. define
\begin{align*}
X_5:=X_4/(1-aT_4^2),\quad Y_5:=Y_4/(1-aT_4^2),\quad T_6:=T_5/(1-aT_4^2).
\end{align*}
Then $K(X_4,Y_4,T_4,T_5)=K(X_5,Y_5,T_4,T_6)$ and the relation (\ref{relYT}) becomes
\begin{align}
3X_5^2-3aY_5^2=1+4T_6+4aT_4^2T_6^2. \label{relZT}
\end{align}

Thus we get $K(x_1,x_2,x_3)^{S_3}=K(X_5,Y_5,T_4T_6,T_6)=K(X_5,Y_5,T_4T_6)$
is rational over $K$ because (\ref{relZT}) becomes linear in $T_6$. \qed

\bigskip
%%%%%%%%%%%%%%%%%%%%%%%%%%%%%%%%%%%%%%%%%%%%%%%%%%%%%%%%%%%%%%%%%%%%%%%%%%%%%%%%%%%%%%%%%%%%%
Proof of Theorem \ref{thmain} when $n=4$ and char $K\neq 2$ ------------------------------------

\medskip

Step $1$. By Lemma \ref{lemactt} we find that
$K(x_1,x_2,x_3,x_4)^{S_4}=K(\alpha)(t_2,t_3,t_4,u)^{\langle\rho\rangle}$
where
\begin{align*}
\rho\ :\ &\alpha\ \longmapsto\ -\alpha,\quad t_2\ \longmapsto\ t_2\, t_3^{-2}\, t_4,\quad
t_3\ \longmapsto\ t_3^{-3}\, t_4^2,\quad t_4\ \longmapsto\ t_3^{-4}\, t_4^3,\\
&u\ \longmapsto\ t_3^{-5}\, t_4^2\, \Delta^2\cdot u^{-1}
\end{align*}
where
\begin{align*}
\Delta^2&=\prod_{1\leq i<j\leq 4}(y_i-y_j)^2\\
&=\ \ t_2^2t_3^2-4t_2^3t_3^2-4t_3^3+18t_2t_3^3-27t_3^4-4t_2^3t_4+16t_2^4t_4
+18t_2t_3t_4-80t_2^2t_3t_4\\
&\hspace*{7mm}-6t_3^2t_4+144t_2t_3^2t_4-27t_4^2+144t_2t_4^2-128t_2^2t_4^2
-192t_3t_4^2+256t_4^3.
\end{align*}

Define $u_1:=t_2$, $u_2:=t_4/t_3$, $u_3:=\rho(t_2)=t_2t_4/t_3^2$.
Then $K(\alpha)(t_2,t_3,t_4,u)$ $=$
$K(\alpha)(u_1,u_2,u_3,u)$ and
\begin{align*}
\rho\ :\ &\alpha\ \longmapsto\ -\alpha,\quad u_1\ \longmapsto\ u_3\ \longmapsto\ u_1,\quad
u_2\ \longmapsto\ u_2,\\
&u\ \longmapsto\ g(u_1,u_2,u_3)\cdot u^{-1}
\end{align*}
where
\begin{align*}
g(u_1,u_2,u_3)=&\frac{u_2}{u_1u_3}
\Bigl(-27u_1^2u_2^2-4u_1u_2u_3+18u_1^2u_2u_3-6u_1u_2^2u_3+144u_1^2u_2^2u_3\\
&-192u_1u_2^3u_3+256u_1u_2^4u_3+u_1^2u_3^2-4u_1^3u_3^2+18u_1u_2u_3^2\\
&-80u_1^2u_2u_3^2-27u_2^2u_3^2+144u_1u_2^2u_3^2-128u_1^2u_2^2u_3^2-4u_1^2u_3^3
+16u_1^3u_3^3\Bigr).
\end{align*}

Define $v_1:=(u_1+u_3)/2$, $v_2:=u_2$, $v_3=\alpha(u_1-u_3)/2$.
Then we obtain $K(\alpha)(u_1,u_2,u_3,u)=K(\alpha)(v_1,v_2,v_3,u)$ and
\begin{align*}
\rho\ :\ &\alpha\ \longmapsto\ -\alpha,\quad v_1\ \longmapsto\ v_1,\quad
v_2\ \longmapsto\ v_2,\quad v_3\ \longmapsto\ v_3,\\
&u\ \longmapsto\ h(v_1,v_2,v_3)\cdot u^{-1}
\end{align*}
where $h(v_1,v_2,v_3)=g(u_1,u_2,u_3)\in K(v_1,v_2,v_3)$ is given as
\begin{align*}
h(v_1,v_2,v_3)=&\frac{v_2}{av_1^2-v_3^2}
\Bigl(av_1^2v_2(-1+4v_1-8v_2)^2(v_1^2-4v_2+4v_1v_2+4v_2^2)\\
&-2v_2v_3^2\bigl(v_1^2-8v_1^3+24v_1^4-2v_2+18v_1v_2-80v_1^2v_2\\
&+24v_2^2+144v_1v_2^2-128v_1^2v_2^2-96v_2^3+128v_2^4\bigr)\\
&-(1/a)v_2v_3^4(-1+8v_1-48v_1^2+80v_2+128v_2^2)-(16/a^2)v_2v_3^6\Bigr).
\end{align*}

%%%%%%%%%%%%%%%%%%%%%%%%%%%%%%%%%%%%%%%%%%%%%%%%%%%%%%%%%%%%%%%%%%%%%%%%%%%%%%%%%%%%%%%%%%%%%
Step $2$. Because $h(v_1,v_2,v_3)$ is still complicated, we define
$p,q,r$ as follows,
\begin{align*}
p:=\frac{1}{2}\Bigl(\frac{1}{u_1}+\frac{1}{u_3}\Bigr)u_2,\quad
q:=\frac{\alpha}{2}\Bigl(\frac{1}{u_1}-\frac{1}{u_3}\Bigr)u_2,\quad
r:=4u_2.
\end{align*}
Then $K(\alpha)(v_1,v_2,v_3,u)=K(\alpha)(p,q,r,u)$. Indeed we have
\begin{align*}
p&=\frac{av_1v_2}{av_1^2-v_3^2},\quad q=-\frac{av_2v_3}{av_1^2-v_3^2},\quad r=4v_2,\\
v_1&=\frac{apr}{4(ap^2-q^2)},\quad v_2=\frac{r}{4},\quad v_3=-\frac{apr}{4(ap^2-q^2)}.
\end{align*}
Hence we obtain $K(x_1,x_2,x_3,x_4)^{S_4}=K(\alpha)(p,q,r,u)^{\langle\rho\rangle}$ and
\begin{align*}
\rho\ :\ &\alpha\ \longmapsto\ -\alpha,\quad p\ \longmapsto\ p,\quad
q\ \longmapsto\ q,\quad r\ \longmapsto\ r,\\
&u\ \longmapsto\ \frac{r^2}{64(ap^2-q^2)^2}\cdot\frac{H(p,q,r)}{u}
\end{align*}
where
\begin{align}
H(p,q,r)={}&a^2(p-r+2pr)^2(-16p^2+r+4pr+4p^2r)\label{eqH}\\
&-a\bigl(-32p^2+r+36pr-12p^2r-20r^2+72pr^2\nonumber\\
&-96p^2r^2-8r^3+32p^2r^3\bigr)q^2+16(-1+r)^3q^4.\nonumber
\end{align}

Define $U:=u\cdot r/(8(ap^2-q^2))$.
Then we get $K(\alpha)(p,q,r,u)=K(\alpha)(p,q,r,U)$ and $\rho$ acts on
$K(\alpha)(p,q,r,U)$ by
\begin{align*}
\rho\ :\ \alpha\ \longmapsto\ -\alpha,\quad p\ \longmapsto\ p,\quad
q\ \longmapsto\ q,\quad r\ \longmapsto\ r,\quad U\ \longmapsto\ H(p,q,r)/U.
\end{align*}
Hence $K(x_1,\ldots,x_4)^{S_4}$ $=$
$K(\alpha)(p,q,r,U)^{\langle\rho\rangle}$ $=$ $K(X,Y,p,q,r)$ where
\[
X=\frac{1}{2}\Bigl(U+\frac{g(p,q,r)}{U}\Bigr),\quad
Y=\frac{1}{2\alpha}\Bigl(U-\frac{g(p,q,r)}{U}\Bigr).
\]
Note that $X$ and $Y$ satisfy the relation
\begin{align}
X^2-aY^2=H(p,q,r).\label{relXYpqr}
\end{align}

%%%%%%%%%%%%%%%%%%%%%%%%%%%%%%%%%%%%%%%%%%%%%%%%%%%%%%%%%%%%%%%%%%%%%%%%%%%%%%%%%%%%%%%%%%%%%
Step $3$.
Because $H(p,q,r)$ in (\ref{eqH}) is biquadratic equation with respect to $q$ and
its constant term has the square factor $(p-r+2pr)^2$, we define
\[
p_2:=p-r+2pr.
\]
Then $p=(p_2+r)/(1+2r)$.
We also define $X_2:=X(1+2r)$, $Y_2:=Y(1+2r)$.
Then $K(x_1,x_2,x_3,x_4)^{S_4}=K(X_2,Y_2,p_2,q,r)$
and the relation (\ref{relXYpqr}) becomes
\begin{align}
X_2^2-aY_2^2={}&a^2p_2^2(-16p_2^2+r-28p_2r+4p_2^2r-8r^2+16p_2r^2+16r^3)\nonumber\\
&-a\bigl(-32p_2^2+r-28p_2r-12p_2^2r-12r^2+120p_2r^2\label{eqX2}\\
&-96p_2^2r^2+48r^3-48p_2r^3+32p_2^2r^3-64r^4+64p_2r^4\bigr)q^2\nonumber\\
&+16(-1+r)^3(1+2r)^2q^4. \nonumber
\end{align}

Note that the right hand side of (\ref{eqX2}) is biquadratic in $q$
with constant term $a^2p_2^2(-16p_2^2+r-28p_2r+4p_2^2r-8r^2+16p_2r^2+16r^3)$.

Hence we define $p_3:=p_2/q$, $X_3:=X_2/q$, $Y_3:=Y_2/q$.
Then the relation (\ref{eqX2}) becomes quadratic in $q$ as follows:
\begin{align}
X_3^2-aY_3^2={}&
ar(-1+4r)^2(-1+ap_3^2+4r)\nonumber\\
&+4ap_3r(7-7ap_3^2-30r+4ap_3^2r+12r^2-16r^3)q\label{relXYpqr2}\\
&+4(-1+ap_3^2-4r-4r^2)(4-4ap_3^2-12r+ap_3^2r+12r^2-4r^3)q^2.\nonumber
\end{align}

Define $q_2:=1/q$, $r_2:=4r$, $X_4:=4X_3/q$, $Y_4:=4Y_3/q$.
Then (\ref{relXYpqr2}) becomes
\begin{align}
X_4^2-aY_4^2={}&4ar_2(-1+r_2)^2(-1+ap_3^2+r_2)q_2^2\nonumber\\
&+4ap_3r_2(28-28ap_3^2-30r_2+4ap_3^2r_2+3r_2^2-r_2^3)q_2\label{relXYpqr3}\\
&(-4+4ap_3^2-4r_2-r_2^2)(64-64ap_3^2-48r_2+4ap_3^2r_2+12r_2^2-r_2^3).\nonumber
\end{align}

Because the relation (\ref{relXYpqr3}) is quadratic with respect to $q_2$, we may
eliminate a linear term of $q_2$ in the usual manner by putting
\[
q_3:=2q_2+\frac{p_3(28-28ap_3^2-30r_2+4ap_3^2r_2+3r_2^2-r_2^3)}
{(-1+r_2)^2(-1+ap_3^2+r_2)}.
\]
Define
\begin{align*}
X_5&:=X_4(-1+r_2)(-1+ap_3^2+r_2),\\
Y_5&:=Y_4(-1+r_2)(-1+ap_3^2+r_2).
\end{align*}
Then the relation (\ref{relXYpqr3}) becomes
\begin{align}
X_5^2-aY_5^2={}&(2+r_2)^2(-1+ap_3^2+r_2)(4-4ap_3^2-5r_2+r_2^2)^3\label{relXYpqr4}\\
&+ar_2(-1+r_2)^4(-1+ap_3^2+r_2)^3q_3^2. \nonumber
\end{align}

Define
\begin{align*}
q_4&:=\frac{q_3(-1+r_2)^2(-1+ap_3^2+r_2)}{(2+r_2)(4-4ap_3^2-5r_2+r_2^2)}
\end{align*}
and
\begin{align*}
X_6&:=\frac{X_5}{(2+r_2)(4-4ap_3^2-5r_2+r_2^2)},\ \
Y_6:=\frac{Y_5}{(2+r_2)(4-4ap_3^2-5r_2+r_2^2)}.
\end{align*}
Then we get $K(x_1,\ldots,x_4)^{S_4}=K(X_6,Y_6,p_3,q_4,r_2)$ and
(\ref{relXYpqr4}) becomes
\begin{align}
X_6^2-aY_6^2=(-1+ap_3^2+r_2)\Bigl((4-4ap_3^2-5r_2+r_2^2)+ar_2q_4^2\Bigr).\label{relXYpqr5}
\end{align}

\bigskip
%%%%%%%%%%%%%%%%%%%%%%%%%%%%%%%%%%%%%%%%%%%%%%%%%%%%%%%%%%%%%%%%%%%%%%%%%%%%%%%%%%%%%%%%%%%%%
Step $4$.
We find the singularities of (\ref{relXYpqr5}).
We get $p_3\pm(1/\sqrt{a})=r_2=X_6=Y_6=0$.
Blow-up with respect to $-1+ap_3^2$, i.e. define
\[
r_3:=r_2/(-1+ap_3^2),\quad X_7:=X_6/(-1+ap_3^2),\quad Y_7:=Y_6/(-1+ap_3^2).
\]
Then $K(p_3,q_4,r_2,X_6,Y_6)=K(p_3,q_4,r_3,X_7,Y_7)$ and
the relation (\ref{relXYpqr5}) becomes
\begin{align}
X_7^2-aY_7^2=(1 + r_3)(-4 - 5r_3 + aq_4^2r_3 - r_3^2 + ap_3^2r_3^2).\label{relXYpqr6}
\end{align}

Define $p_4:=p_3r_3$.
Then the relation (\ref{relXYpqr6}) becomes
\begin{align}
X_7^2-aY_7^2=(1 + r_3)(-4 - 5r_3  + aq_4^2r_3 - r_3^2+ ap_4^2).\label{relXYpqr7}
\end{align}

\bigskip
%%%%%%%%%%%%%%%%%%%%%%%%%%%%%%%%%%%%%%%%%%%%%%%%%%%%%%%%%%%%%%%%%%%%%%%%%%%%%%%%%%%%%%%%%%%%%
Step $5$.
The equation (\ref{relXYpqr7}) still has a singular locus $p_4\pm q_4=r_3+1=X_7=Y_7=0$.
Define $p_5:=p_4+q_4$, $r_4:=r_3+1$.
Then the relation (\ref{relXYpqr7}) becomes
\begin{align}
X_7^2-aY_7^2=r_4(ap_5^2 - 2ap_5q_4 - 3r_4 + aq_4^2r_4 - r_4^2)\label{relXYpqr8}
\end{align}
with singular locus $S=(p_5=r_4=X_7=Y_7=0)$.
Blow-up (\ref{relXYpqr8}) along $S$, i.e. define
$r_5:=r_4/p_5$, $X_8:=X_7/p_5$, $Y_8:=Y_7/p_5$.
Then the relation (\ref{relXYpqr8}) becomes
\[
X_8^2-aY_8^2=r_5(ap_5 - 2aq_4 - 3r_5 + aq_4^2r_5 - p_5r_5^2).
\]
Note that the relation above is linear in $p_5$.
Hence we conclude that the fixed field $K(x_1,\ldots,x_4)^{S_4}$ $=$
$K(X_8,Y_8,q_4,r_5)$ is rational over $K$. \qed

\bigskip

%%%%%%%%%%%%%%%%%%%%%%%%%%%%%%%%%%%%%%%%%%%%%%%%%%%%%%%%%%%%%%%%%%%%%%%%%%%%%%%%%%%%%%%%%%%%%
\section{Using the structures of rings of invariants}\label{se4}

In this section, we give an another proof of Theorem \ref{thmain}
in the case of $n=3, 4$ and char $K=2$ by using the structure of
$K(x_1, \cdots, x_n)^{A_n}$.

Throughout this section, we assume that char $K=2$.

Let $s_i$ be the elementary symmetric function in $x_1,\ldots,x_n$
of degree $i$, $1\leq i\leq n$.

Since  the invariant ring $K[x_1,\ldots,x_n]^{A_n}$ is a free
module of rank $2$ over the subring $K[x_1,\ldots,x_n]^{S_n}$ $=$
$K[s_1,\ldots,s_n]$ by Revoy's Theorem \cite{Re} (see also \cite[
Example 1, page 75]{NS}), we will find a free basis of
$K[x_1,\ldots,x_n]^{A_n}$ over $K[x_1,\ldots,x_n]^{S_n}$. For
$n=3$ (resp. $n=4$), it suffices to find a polynomial of degree
$3$ (resp. $6$), which is in $K[x_1,\ldots,x_n]^{A_n}$, but
doesn't belong to $K[x_1,\ldots,x_n]^{S_n}$ by \cite[ Example 1,
page 75]{NS}.

Define
\begin{align*}
b_3\, &:=\,\sum_{\sigma\in A_3} \sigma (x_1 x_2^2)
=x_1 x_2^2+x_2 x_3^2+x_3 x_1^2,\\
b_4\, &:=\,\sum_{\sigma\in A_4} \sigma (x_1 x_2^2 x_3^3)\\
&\ =\, x_1^2 x_2^3 x_3+x_1^3 x_2 x_3^2+x_1 x_2^2 x_3^3+x_1^3 x_2^2
x_4+x_2^3 x_3^2 x_4
+x_1^2 x_3^3 x_4\\
&\ \ \quad +x_1 x_2^3 x_4^2+x_1^3 x_3 x_4^2+x_2 x_3^3 x_4^2+x_1^2
x_2 x_4^3 +x_2^2 x_3 x_4^3+x_1 x_3^2 x_4^3.
\end{align*}

For $n=3$ (resp. $n=4$), it follows that $\{1,b_3\}$ (resp.
$\{1,b_4\}$) is a free basis of $K[x_1,x_2,x_3]^{A_3}$ (resp.
$K[x_1,x_2,x_3,x_4]^{A_4}$), i.e.
\begin{align*}
K[x_1,x_2,x_3]^{A_3}&=K[s_1,s_2,s_3]\oplus b_3\,
K[s_1,s_2,s_3],\\
K[x_1,x_2,x_3,x_4]^{A_4}&=K[s_1,s_2,s_3,s_4]\oplus b_4\,
K[s_1,s_2,s_3,s_4].
\end{align*}

We conclude that

\begin{lemma}
Let $K$ be a field of char $K=2$ and write $S_3=\langle
A_3,\tau\rangle$, $S_4= \langle A_4, \tau\rangle$
where $\tau=(12)$.\\
{\rm (i)} For $n=3$, we have
\[
K(x_1,x_2,x_3)^{A_3}=K(s_1,s_2,s_3,b_3)
\]
and $s_1,s_2,s_3,b_3$ satisfy the following relation:
\[
b_3^2 + b_3\, s_1s_2 + s_2^3 + b_3\, s_3 + s_1^3 s_3 + s_3^2=0.
\]

{\rm (ii)} For $n=4$, we have
\[
K(x_1,x_2,x_3,x_4)^{A_4}=K(s_1,s_2,s_3,s_4,b_4)
\]
and $s_1,s_2,s_3,s_4,b_4$ satisfy the following relation:
\[
b_4^2 + b_4\, s_1 s_2 s_3 + b_4\, s_3^2 + s_2^3 s_3^2 + s_1^3
s_3^3 + s_3^4 + b_4\, s_1^2 s_4 + s_1^2 s_2^3 s_4 + s_1^4 s_4^2=0.
\]
\end{lemma}

\bigskip
%%%%%%%%%%%%%%%%%%%%%%%%%%%%%%%%%%%%%%%%%%%%%%%%%%%%%%%%%%%%%%%%%%%%%%%%%%%%%%%%%%%%%%%%%%%%%
Proof of Theorem \ref{thmain} when $n=3$ and char $K=2$
------------------------

Note that $\tau$ acts on $K(x_1,x_2,x_3)^{A_3}=K(s_1,s_2,s_3,b_3)$
as follows
\begin{align*}
&s_1\longmapsto \frac{a s_2}{s_3},\ \ s_2\longmapsto
\frac{a^2s_1}{s_3},\ \ s_3\longmapsto \frac{a^3}{s_3}, \ \
b_3\longmapsto\frac{a^3b_3}{s_3^2}.
\end{align*}

Apply Theorem \ref{lemMasuda}. We find
$K(x_1,x_2,x_3)^{A_3}=K(s_3,u,v)$ where $u$ and $v$ are the same
as in Theorem \ref{lemMasuda}. It is not difficult to verify that
\begin{align*}
u=\frac{b_3+s_3}{s_1^2+s_2},\quad v=\frac{b_3+s_1s_2}{s_1^2+s_2}.
\end{align*}

Moreover, the action of $\tau$ is given by
\begin{align*}
\tau\ :\ &\ s_3\ \longmapsto\ \frac{a^3}{s_3},\quad u\
\longmapsto\ \frac{au}{u^2-uv+v^2},\quad v\ \longmapsto\
\frac{av}{u^2-uv+v^2}.
\end{align*}

Define $w:=u/v$. Then $K(x_1,x_2,x_3)^{A_3}=K(s_3,v,w)$ and
\begin{align*}
\tau\ :\ &\ s_3\ \longmapsto\ \frac{a^3}{s_3},\quad v\
\longmapsto\ \frac{a}{v(1-w+w^2)},\quad w\ \longmapsto\ w.
\end{align*}
By Theorem \ref{lemHK}, $K(x_1,x_2,x_3)^{S_3}$ $=$
$K(s_3,v,w)^{\langle\tau\rangle}$ is rational over $K$.\qed

\bigskip
%%%%%%%%%%%%%%%%%%%%%%%%%%%%%%%%%%%%%%%%%%%%%%%%%%%%%%%%%%%%%%%%%%%%%%%%%%%%%%%%%%%%%%%%%%%%%
Proof of Theorem \ref{thmain} when $n=4$ and char $K=2$
------------------------

For $n=4$, $\tau$ acts on
$K(x_1,x_2,x_3,x_4)^{A_4}=K(s_1,s_2,s_3,s_4,b_4)$ as
\begin{align*}
&s_1\longmapsto \frac{a s_3}{s_4},\ \ s_2\longmapsto \frac{a^2
s_2}{s_4},\ \
s_3\longmapsto \frac{a^3 s_1}{s_4},\\
&s_4\longmapsto \frac{a^4}{s_4},\ \ b_4\longmapsto
\frac{a^6(b_4+s_1s_2s_3+s_3^2+s_1^2s_4)}{s_4^3}.
\end{align*}

Define
\begin{align*}
t_1:=\frac{s_1s_3}{s_2},\quad t_2:=s_2,\quad t_3=s_3,\quad
t_4:=\frac{s_1s_2s_3+s_3^2+s_1^2s_4}{s_2^2},\quad
t_5:=\frac{b_4+s_2^3}{s_2}.
\end{align*}

It follows that $K(s_1,s_2,s_3,s_4,b_4)=K(t_1,t_2,t_3,t_4,t_5)$.
It is easy to check that the relation among the generators
$t_1,\ldots,t_5$ is given by
\begin{align*}
t_1^3 + t_1^2t_2 + t_1t_2^2 + t_2^3 + t_2t_4^2 + t_2t_4t_5 +
t_2t_5^2=0.
\end{align*}
Define
\[
u_1:=t_1,\quad u_2:=t_2/t_1,\quad u_3:=t_3, \quad
u_4:=t_4/(t_1+t_2),\quad u_5:=t_5/(t_1+t_2).
\]
Then we get $K(t_1,\ldots,t_5)=K(u_1,\ldots,u_5)$ with the
relation
\[
u_2(u_4^2+u_4u_5+u_5^2+1)+1=0.
\]
Because this relation is linear in $u_2$, we obtain the following
lemma.

\begin{lemma}
$K(x_1,\ldots,x_4)^{A_4}=K(u_1,u_3,u_4,u_5)$ where
\begin{align*}
u_1=\frac{s_1s_3}{s_2},\quad u_3=s_3,\quad
u_4=\frac{s_1s_2s_3+s_3^2+s_1^2s_4}{s_2(s_2^2+s_1s_3)},\quad
u_5=\frac{b_4+s_2^3}{s_2(s_2^2+s_1s_3)}.
\end{align*}
\end{lemma}

\bigskip

Now we will prove Theorem \ref{thmain} when $n=4$ and char $K=2$.

\medskip

Write $p=u_1, q=u_3, r=u_4, s=u_5$ and $\tau=(12)\in S_4\setminus
A_4$. Note that
$K(x_1,\ldots,x_4)^{S_4}=K(p,q,r,s)^{\langle\tau\rangle}$ and the
action of $\tau$ on $K(p,q,r,s)$ is given by
\begin{align*}
\tau\ :\ p&\ \longmapsto\ \frac{r^2+rs+s^2+1}{ap},\\
q&\ \longmapsto
\frac{a^3p^6q}{(r^2+rs+s^2+1)^3+p^3q\bigl((r+1)(r^2+rs+s^2+1)+1\bigr)},\\
r&\ \longmapsto r,\quad s \longmapsto s+r.
\end{align*}
Define
\[
t:=\frac{(r^2+rs+s^2+1)^3}{p^3q\bigl((r+1)(r^2+rs+s^2+1)+1\bigr)}.
\]
Then we have
$K(x_1,x_2,x_3,x_4)^{S_4}=K(p,q,r,s)^{\langle\tau\rangle}$ $=$
$K(p,t,r,s)^{\langle\tau\rangle}$ and the action of $\tau$ on
$K(p,t,r,s)$ is given by
\begin{align*}
\tau\ :\ p\ \longmapsto\ \frac{r^2+rs+s^2+1}{ap},\quad t\
\longmapsto\ t+1,\quad r\ \longmapsto\ r,\quad s\ \longmapsto\
s+r.
\end{align*}
Define
\[
A:=r+s+rt,\quad B:=(r+s)/s,\quad C:=pr/s.
\]
It follows that $K(p,q,r,s)=K(r,A,B,C)$. Thus
$K(x_1,x_2,x_3,x_4)^{S_4}=K(r)(A$, \ $B,C)^{\langle\tau\rangle}$
and
\[
\tau\ :\ r\ \longmapsto\ r,\quad A\ \longmapsto\ A,\quad B\
\longmapsto\ \frac{1}{B},\quad C\ \longmapsto\
\frac{1}{a}\Bigl((r^2+1)\bigl(\frac{1}{B}+B\bigl)+r^2\Bigr)\Big/C.
\]
Apply Theorem \ref{lemHK}. We find that $K(x_1,x_2,x_3,x_4)^{S_4}$
is rational over $K$. \qed

\newpage

%%%%%%%%%%%%%%%%%%%%%%%%%%%%%%%%%%%%%%%%%%%%%%%%%%%%%%%%%%%%%%%%%%%%%%%%%%%%%%%%%%%%%%%%%%%%%

\end{document}